\documentclass[reqno]{amsart}
\usepackage{amsmath,amsthm,amssymb,graphicx}
\usepackage{latexsym}
\usepackage{eucal}

\newtheorem{theorem}{Theorem}[section]
\newtheorem{lemma}{Lemma}[section]

\newtheorem{conjecture}{Conjecture}[section]

\theoremstyle{definition}

\let\Re=\undefined
\DeclareMathOperator{\Re}{Re}
\let\Im=\undefined
\DeclareMathOperator{\Im}{Im}

\begin{document}
\title[On a conjecture by Y. Last \ldots]{On a conjecture by Y. Last}
\author{Sergey A. Denisov}

\address{University of Wisconsin-Madison,
Mathematics Department, 480 Lincoln Dr., \linebreak Madison, WI,
53706-1388, USA,  {\rm e-mail: denissov@math.wisc.edu}}
%\address{Mathematics, 253-37, Caltech, Pasadena, CA, 91125, USA,
%{\rm e-mail: denissov@caltech.edu}}
\maketitle \centerline{\it{To the memory of George G. Lorentz}}
\begin{abstract}
We prove a conjecture due to Y. Last. The new determinantal
representation for transmission coefficient of Jacobi matrix is
obtained.
\end{abstract} \vspace{1cm}

In this paper we consider discrete Schr\"odinger operator on the
half-lattice with bounded real potential $v$

\begin{equation}
 {J}=\left[
\begin{array}{ccccc}
 {v}_{1} &  1 & 0 & 0 & \ldots  \\
 1 &  {v}_{2} &  1 & 0 & \ldots  \\
0 &  1 &  {v}_{3} &  1 & \ldots  \\
0 & 0 &  1 &  {v}_{4} & \ldots  \\
\ldots  & \ldots  & \ldots  & \ldots  & \ldots
\end{array}%
\right] \label{jacobi}
\end{equation}

In \cite{last}, Y. Last posed the following problem
\begin{conjecture}
Do the following conditions: $v_n\to 0$ and
\[
v_{n+q}-v_n\in \ell^2
\]
($q\in \mathbb{Z}^+$-- fixed) guarantee that
$\sigma_{ac}(J)=[-2,2]$?
\end{conjecture}
The symbol $\sigma_{ac}(J)$ conventionally denotes the absolutely
continuous (a.c.) spectrum of self-adjoint operator $J$. In this
paper, we give an affirmative answer to this question. The
manuscript consists of two sections. The first one is mostly
algebraic, it contains the determinantal formula for the so-called
transmission coefficient that allows us to immediately treat the
case $q=1$. In the second part, we show how asymptotical methods
for difference equations provide the solution for any $q$. The
appendix contains an
 elementary lemmas for harmonic functions which are
used in the paper.

 Recently, many results
on the characterization of parameters in the Jacobi matrix through
the spectral data were obtained and the $\ell^2$--condition on
coefficients was often involved in one form or another (see, e.g.,
\cite{DK, ks, k, zs}). This paper makes the next step in this
direction by developing the technique suggested in \cite{d}.

We will use notations: $(\delta v)_n=v_{n+1}-v_n,
(\delta^{(q)}v)_n=v_{n+q}-v_n$, $\chi_{j\in M}$ is the
characteristic function of the set $M$. For the sequence $\alpha\in
\ell^p$, the symbol $\|\alpha\|_p$ denotes its norm in $\ell^p$. As
usual, the symbol $C$ denotes the positive constant which might take
different values in different formulas. For any matrix $B\in
\mathbb{C}^{k\times k}$, the symbol $\|B\|$ will denote its operator
norm in $\mathbb{C}^k$. Consider a linear bounded operator $A$
acting in the Hilbert space. Assume that it is Hilbert-Schmidt, i.e.
$A\in S_2$. Then,  we define the regularized determinant by the
formula (see, e.g., \cite{sim})
\[
\det{}_2 (I+A)=\det (I+R_2(A))
\]
where
\[R_2(A)=(I+A)\exp(-A)-I\in S_1
\]
The symbol $S_p$ is reserved for the standard Schatten--von Neumann
class.
\section{Determinantal formula  and $q=1$.}

It will be convenient for us to start with Jacobi matrices on
$\ell^2(\mathbb{Z})$.  Let $H$ be a discrete Schr\"odinger operator
on $\ell^2(\mathbb{Z})$ with potential $v$.  Later on, we will make
the following choice for $v$. For positive indices, it will be taken
from (\ref{jacobi}). For negative indices, it will be set to zero.
By $H_0$ we will denote the ``free" case, i.e. the case when the
potential $v$ is identically zero on all of $\mathbb{Z}$. Consider
the right and the left shifts acting on $\ell^2(\mathbb{Z})$
\[
(Rf)_n=f_{n-1}, \quad L=R^*
\]
Obviously, $R=L^{-1}$. For any $z\in \mathbb{C}$, introduce the
following diagonal operators $\Lambda(z)=\{\lambda_n(z)\}$ and
$\Lambda^0(z)=\tilde\lambda(z)\cdot I$, where
\begin{equation}
\lambda_n(z)=\frac{z-v_n-\left[(z-v_n)^2-4\right]^{1/2}}{2},\quad
\tilde\lambda(z)=\frac{z-\left[z^2-4\right]^{1/2}}{2}
\label{lam-poz}
\end{equation}
and $\sqrt{\,}$ has a cut along the positive axis. Notice that the
both functions map the upper half--plane into the lower half of the
unit disc, i.e. $\{z\in \mathbb{C}: |z|<1, \Im z<0\}$. Let
$\delta\Lambda(z)=\Lambda(z)-\Lambda_0(z)$. We also need the finite
dimensional versions of these operators. Take $m=2n+1$--dimensional
linear space span$\{e_{-n},\ldots, e_n\}$ and let $R_n$ and $L_n$ be
right and left cyclic shifts, respectively. They are unitary
operators and $L_n=R_n^{-1}$. Let $\Lambda_n$, $\Lambda^0_n$,
$\delta \Lambda_n=\Lambda_n-\Lambda^0_n$ be restrictions of
$\Lambda$, $\Lambda^0$, and $\delta \Lambda$. Define
$\omega_j=\lambda_{j+1}-\lambda_j, j=-n,\ldots, n-1$,
$\omega_n=\lambda_{-n}-\lambda_n$ and the corresponding diagonal
operator $\Omega_n=\{\omega_j\}, j=-n, \ldots, n$.

 If $K_n=L_n+R_n-\Lambda_n-\Lambda_n^{-1}$, then we have an elementary

\begin{lemma}
For any $z\notin \mathbb{R}$, we have
\[
(L_n-\Lambda_n)(R_n-\Lambda_n)=-\Lambda_n
K_n-\Omega_nL_n=-K_n\Lambda_n-R_n\Omega_n
\]
\[
(R_n-\Lambda_n)(L_n-\Lambda_n)=-K_n\Lambda_n+\Omega_nL_n
\]
\end{lemma}
\begin{proof}
The proof is a straightforward calculation.
\end{proof}
\begin{lemma}
For any $z\in \mathbb{C}^+$,
\begin{equation}
\begin{array}{l}
\displaystyle \det
K_n=

-\frac{\det(L_n-\Lambda_n)\det(R_n-\Lambda_n)}{\det\Lambda_n}\\
\displaystyle \times \exp \Bigl({\rm tr}\left[
\Omega_nL_n[(L_n-\Lambda_n)(R_n-\Lambda_n)]^{-1}\right]\Bigr)\\
\displaystyle \times\det{}_2 [I+\Omega_nL_n
[(L_n-\Lambda_n)(R_n-\Lambda_n)]^{-1}]\\
\\

=\displaystyle -\frac{\det(L_n-\Lambda_n)\det(R_n-\Lambda_n)}{\det\Lambda_n}\\
\displaystyle \times \exp \Bigl({\rm tr}\left[
\Omega_nL_n[(R_n-\Lambda_n)(L_n-\Lambda_n)-R_n\Omega_n-\Omega_nL_n]^{-1}\right]\Bigr)\\
\displaystyle \times\det{}_2 [I+\Omega_nL_n
[(L_n-\Lambda_n)(R_n-\Lambda_n)]^{-1}]
\end{array}\label{j1}
\end{equation}
\begin{equation}
\begin{array}{l}
\displaystyle
 \det
K_n=-\frac{\det(L_n-\Lambda_n)\det(R_n-\Lambda_n)}{\det\Lambda_n}\\
\times \exp \Bigl({\rm tr}\left[
-\Omega_nL_n[(R_n-\Lambda_n)(L_n-\Lambda_n)]^{-1}\right]\Bigr)
\times\det{}_2 [I-\Omega_nL_n
[(R_n-\Lambda_n)(L_n-\Lambda_n)]^{-1}]
\end{array}\label{j2}
\end{equation}

\begin{equation}
\begin{array}{l}
\displaystyle \det
K_n=-\frac{\det(L_n-\Lambda_n)\det(R_n-\Lambda_n)}{\det\Lambda_n}\\
\displaystyle  \times \exp \left(\frac 12 \,{\rm tr}\left[
\Omega_nL_n(L_n-\Lambda_n)^{-1}(R_n-\Lambda_n)^{-1}(R_n\Omega_n+\Omega_nL_n)
(R_n-\Lambda_n)^{-1}(L_n-\Lambda_n)^{-1}\right]\right)\\
\times \Bigl[ \det{}_2 [I+\Omega_nL_n
(R_n-\Lambda_n)^{-1}(L_n-\Lambda_n)^{-1}]\det{}_2 [I-\Omega_nL_n
(L_n-\Lambda_n)^{-1}(R_n-\Lambda_n)^{-1}\Bigr]^{1/2}
\end{array}\label{j3}
\end{equation}
\end{lemma}

\begin{proof}
(\ref{j1}) and (\ref{j2}) follow immediately from the previous
lemma. Multiplication of (\ref{j1}) and (\ref{j2}) yields
(\ref{j3}) after application of the second resolvent identity:
\[
(A+V)^{-1}-A^{-1}=-A^{-1}V(A+V)^{-1}
\]
and taking the suitable square root.
\end{proof}

Writing down the formula (\ref{j3}) for the ``free" case with
\[
K_n^0=L_n+R_n-\Lambda_n^0-(\Lambda_n^0)^{-1}
\]
and dividing, we have
\begin{equation}
\begin{array}{l}
\displaystyle \det
\left[K_n/K_n^0\right]=\frac{\det(I-(L_n-\Lambda_n^0)^{-1}\delta\Lambda_n)
\det(I-(R_n-\Lambda_n^0)^{-1}\delta\Lambda_n)}{\det\left[\Lambda_n/\Lambda_n^0\right]}\\
\displaystyle  \times \exp \left(\frac 12 \,{\rm tr}\left[
\Omega_nL_n(L_n-\Lambda_n)^{-1}(R_n-\Lambda_n)^{-1}(R_n\Omega_n+\Omega_nL_n)
(R_n-\Lambda_n)^{-1}(L_n-\Lambda_n)^{-1}\right]\right)\\
\times \Bigl[ \det{}_2 [I+\Omega_nL_n
(R_n-\Lambda_n)^{-1}(L_n-\Lambda_n)^{-1}]\det{}_2 [I-\Omega_nL_n
(L_n-\Lambda_n)^{-1}(R_n-\Lambda_n)^{-1}\Bigr]^{1/2}
\end{array}\label{j4}
\end{equation}
Later on we will need the following bound
\begin{lemma}
For any $z\in \mathbb{C}^+$ and $v$, we have the following estimates
for the operator norms
\begin{equation}
\|(L_n-\Lambda_n)^{-1}\|\leq C(\Im z)^{-1}(1+\Im z),\quad
\|(L-\Lambda)^{-1}\|\leq C(\Im z)^{-1}(1+\Im z) \label{unifb1}
\end{equation}
\[
\|(R_n-\Lambda_n)^{-1}\|\leq C(\Im z)^{-1}(1+\Im z),\quad
\|(R-\Lambda)^{-1}\|\leq C(\Im z)^{-1}(1+\Im z)
\]\label{llls}
\end{lemma}
\begin{proof}
Writing $\tilde{\lambda}(z)$ in polar coordinates, one can prove
that
\[
|\tilde\lambda(z)|=\left|\frac{z-\sqrt{z^2-4}}{2}\right|\leq
\frac{\sqrt{4+(\Im z)^2}-\Im z}{2}, \quad z\in \mathbb{C}^+
\]
and
\[
\|\Lambda_n\|\leq \frac{\sqrt{4+(\Im z)^2}-\Im z}{2}
\]
Since
\[
\frac{C\Im z}{1+\Im z}<1-\|\Lambda_n\|\leq |((L_n-\Lambda_n)f, L_n
f)|\leq \|(L_n-\Lambda_n)f\|, \quad \|f\|=1
\]
we have the statement of the lemma. The statements for $L, R_n, R$
have the same proofs.
\end{proof}

Taking $n\to\infty$ in (\ref{j4}), we get
\begin{lemma}
Assume that $v$ is supported on $|j|\leq l$. For $z\in
 \mathbb{C}^+$,
\begin{equation}
\begin{array}{l}
\displaystyle \det
\left[(H-z)/(H_0-z)\right]=\frac{\det(I-(L-\Lambda^0)^{-1}\delta\Lambda)
\det(I-(R-\Lambda^0)^{-1}\delta\Lambda)}{\det\left[\Lambda/\Lambda^0\right]}\\
\displaystyle  \times \exp \left(\frac 12 \,{\rm tr}\left[ \Omega
L(L-\Lambda)^{-1}(R-\Lambda)^{-1}(R\Omega+\Omega L)
(R-\Lambda)^{-1}(L-\Lambda)^{-1}\right]\right)\\
\times \Bigl[ \det{}_2 [I+\Omega L
(R-\Lambda)^{-1}(L-\Lambda)^{-1}]\det{}_2 [I-\Omega L
(L-\Lambda)^{-1}(R-\Lambda)^{-1}\Bigr]^{1/2}
\end{array}\label{j5}
\end{equation}

\end{lemma}
\begin{proof}
Since $v$ is compactly supported, all determinants and traces in
(\ref{j4}) and (\ref{j5}) are taken of the finite matrices with size
of order $\sim l$. Therefore, it is sufficient to check that the
corresponding matrices converge componentwise, which follows from
the bound (\ref{unifb1}).  In the same way one can show that
\[
\det [I+(K_n^0)^{-1}(K_n-K_n^0)] \to \det [I+(H_0-z)^{-1}(H-H_0)]
\]
Then, taking $n\to\infty$ in (\ref{j4}), one has (\ref{j5}).
\end{proof}

\begin{lemma} If $z\in \mathbb{C}^+$, then
\[
[(L-\Lambda^0)^{-1}f]_n=\sum\limits_{k=-\infty}^n
\tilde\lambda^{n-k}(z)f_{k-1},\quad
[(R-\Lambda^0)^{-1}f]_n=\sum\limits_{k=n}^\infty
\tilde\lambda^{k-n}(z)f_{k+1}
\]
\label{res1}
\end{lemma}
\begin{proof}
The sums are in $\ell^2(\mathbb{Z})$ by Young's inequality for
convolutions since \mbox{$|\tilde\lambda(z)|<1$}. The rest is a
direct calculation.
\end{proof}
By inspection, we have
\begin{lemma}
If $z\in \mathbb{C}^+$, then
\[
\det\left[ I- (L-\Lambda^0)^{-1}\delta \Lambda\right]=1,\quad
\det\left[ I- (R-\Lambda^0)^{-1}\delta \Lambda\right]=1
\]
\[
\det \left[\Lambda^0\Lambda^{-1}\right]=\prod\limits_{j=-l}^l
\left[{\tilde\lambda(z)}/{\lambda_j(z)}\right]
\]
\end{lemma}
\begin{proof}
It is a direct corollary from lemma \ref{res1}.
\end{proof}
Thus the formula (\ref{j5}) can be simplified and we have
\begin{lemma}
For compactly supported $v$,
\begin{equation}
\begin{array}{l}
\displaystyle \det
\left[(H-z)/(H_0-z)\right]=\prod\limits_{j=-\infty}^\infty
\left[{\tilde\lambda(z)}/{\lambda_j(z)}\right]\\\displaystyle \times
\exp \left(\frac 12 \,{\rm tr}\left[ \Omega
L(L-\Lambda)^{-1}(R-\Lambda)^{-1}(R\Omega+\Omega L)
(R-\Lambda)^{-1}(L-\Lambda)^{-1}\right]\right)\\
\times \Bigl[ \det{}_2 [I+\Omega L
(R-\Lambda)^{-1}(L-\Lambda)^{-1}]\det{}_2 [I-\Omega L
(L-\Lambda)^{-1}(R-\Lambda)^{-1}\Bigr]^{1/2}
\end{array}\label{j6}
\end{equation}
\end{lemma}

Consider the first factor. Take any sequence $v=\{v_j\}, j\in
\mathbb{Z}$ of real numbers such that $v_j\to 0$ as $|j|\to \infty$
and let $v^N=v\cdot \chi_{|j|<N}$ be its truncation where $N$ is
large. For each $N$, introduce\footnote{\,The use of symbol $WKB$ is
justified by analogous results in asymptotical theory of ordinary
differential equation (see, e.g., \cite{fed}). This abbreviation is
after Wentzel--Kramers--Brillouin.}

\[
WKB_{v^N}(z)=\prod\limits_{j=-N}^N
\left[{\tilde\lambda(z)}/{\lambda_j(z)}\right]
\]
where $\lambda_j(z)$ and $\tilde\lambda(z)$ are defined in
(\ref{lam-poz}). Obviously, if $v$ is compactly supported, then
$WKB_{v^N}(z)$ will coincide with the first factor in (\ref{j6}) as
long as $N$ is large.

Notice that, for any fixed $\delta$, we have
\begin{equation}
|WKB_{v^N}(z)|\sim 1\label{in1}
\end{equation}
for any $z\in [-2+\delta, 2-\delta]$ and any $N$ provided that
$\|v\|_\infty$ is small.

In scattering theory, the inverse to the transmission coefficient is
usually denoted by $a(z)$. We will introduce its modification (or,
rather, regularization). It will be denoted by $a_m(z)$. For
compactly supported $v$, consider the second and the third factors
in (\ref{j6}). Let

$a_m(z)=$
\begin{equation}
\begin{array}{l} \displaystyle  =\exp \left(\frac 12 \,{\rm
tr}\left[ \Omega L(L-\Lambda)^{-1}(R-\Lambda)^{-1}(R\Omega+\Omega
L)
(R-\Lambda)^{-1}(L-\Lambda)^{-1}\right]\right)\\
\times \Bigl[ \det{}_2 [I+\Omega L
(R-\Lambda)^{-1}(L-\Lambda)^{-1}]\det{}_2 [I-\Omega L
(L-\Lambda)^{-1}(R-\Lambda)^{-1}\Bigr]^{1/2}
\end{array}\label{j7}
\end{equation}
The relevance of $a_m(z)$ to the scattering will be clear from the
proof of the Theorem~\ref{th1}.

We want to control $a_m(z)$ for $\Im z>0$ and $|z|<4$. Specifically,
we need estimates on the boundary behavior as $z$ approaches
$[-2,2]$. For the Hilbert--Schmidt norm $\|\Omega\|_{S^2}$ of
$\Omega$, we have
\[
\|\Omega\|_{S^2} \leq C\|\delta v\|_2
\]
Combining this estimate with (\ref{unifb1}), we get
\begin{lemma}
For $\Im z>0,\, |z|<4$,
\begin{equation}
\ln |a_m(z)|\leq C\frac{\|\delta v\|_2^2}{(\Im z)^4}\label{main}
\end{equation}
Also, for any fixed $\epsilon>0$ and any $z: \Im z>\epsilon, \,
|z|<4$, we have
\begin{equation}
|a_m(z)|>C(\epsilon,\|\delta v\|_2)>0 \label{second}
\end{equation}
provided that $\|\delta v\|_2$ is small.
\end{lemma}
\begin{proof}
The estimates follow from the Lemma \ref{llls} and from the
properties of the trace and $\det_{\,2}$ (see \cite{sim}, p.107 (b)
and \cite{sim4}, p.242 problem 28.2).
\end{proof}
Now, we are ready to prove the main statement of the first
section.
\begin{theorem}
Assume that $v_n\to 0$ and $v_{n+1}-v_n\in \ell^2(\mathbb{Z}^+)$.
Then, $\sigma_{ac}(J)=[-2,2]$.\label{th1}
\end{theorem}
\begin{proof}
By Weyl's theorem, the essential spectrum of $J$ is $[-2,2]$. By
the Kato-Rosenblum theorem, the support of a.c. spectrum does not
change under the trace-class perturbations and
$\sigma_{ac}(J)=\sigma_{ac}(J(L))$ where $J(L)$ has potential
$v(L)=v\cdot \chi_{j>L}$.
\begin{equation}
\|v(L)\|_\infty+\|\delta v(L)\|_2\to 0, \quad L\to\infty
\label{third}
\end{equation}
and therefore we can assume $\|v\|_\infty+\|\delta v\|_2$ to be as
small as we wish.

 For
large $N$, consider truncations $v^N$, i.e. $v^N=v\cdot \chi_{j<N}$.
Also, take $H^N$ on $\ell^2(\mathbb{Z})$ with potential $v_j=v^N_j$
for $j\geq 0$ and $v_j=0$ for $j<0$. Let $J^N$ denotes analogous
truncation for $J$. Then, the Jost function $\psi_n(k)$ is defined
as the solution to $H^N\psi=(k+k^{-1})\psi$ that satisfies
$\psi_n(k)=k^n$ for $n>N$. It is well-known that such a solution
exists for all $k\neq 0$. We will be interested in $k:\Im k\leq 0,
0<|k|\leq 1$ which corresponds to $z=k+k^{-1}\in
\overline{\mathbb{C}^+}$.

Since $v^N_j=0$ for negative $j$, $\psi_n(k)=a^N(z)k^n+b^N(z)k^{-n}$
for $n<0$. We will use the following well-known facts
\[
a^N(z)=\det\left[\frac{H^N-z}{H_0^N-z}\right]
\]
\begin{equation}
|a^N(z)|\geq 1\,\, {\rm for }\,\, z\in (-2,2) \label{posit}
\end{equation}
\begin{equation}
\frac{1}{|a^N(z)|^2}=\frac{4|\sin\theta|}{|m^N(z)+e^{i\theta}|^2}\Im
m^N(z),\quad z=2\cos\theta\in [-2,2]\label{svyaz2}
\end{equation}
(see p. 346, \cite{DK}), here $m^N(z)$ is the Stieltjes transform of
the spectral measure $d\rho_N(\lambda)$ of $J^N$. From
(\ref{posit}),
\[
\int\limits_{-2+\delta}^{2-\delta} \ln |a^N(z)|dz>0
\]
for any small $\delta>0$.

  Consider the
formula (\ref{j6}) for $a^N(z)$ and use (\ref{in1}), (\ref{j7}), and
Lemma \ref{lemma2} from Appendix with $f(z)=-\ln|a_m^N(z)|$,
$a=-2+\delta$, $b=2-\delta$. The estimates (\ref{main}),
(\ref{second}), and (\ref{third}) guarantee its applicability
because the function $a^N(z)$ is continuous up to the real line in
the specified domain since $v^N$ has a finite support. Thus, we have
\[
-\int\limits_{-2+\delta}^{2-\delta} \ln |a^N_m(z)|dz>-C
\]
Therefore, due to (\ref{in1}) and (\ref{svyaz2}),
\begin{equation}
\int\limits_{-2+\delta}^{2-\delta} \ln \rho_N'(z)dz\geq -C
\end{equation}
uniformly in $N$. Since $d\rho_N(\lambda)\to d\rho(\lambda)$ in the
weak sense \cite{DK}, the semicontinuity argument from \cite{ks},
Corollary 5.3 gives
\begin{equation}
\int\limits_{-2+\delta}^{2-\delta} \ln \rho'(z)dz\geq -C
\end{equation}
for all $\delta>0$. That implies that the a.c. part of the spectrum
covers $[-2,2]$.
\end{proof}

\section{Last's conjecture for any $q$.}

In this section, we will apply the standard method of asymptotical
analysis to study the  Schr\"odinger difference relation, then the
asymptotics obtained will be analyzed to conclude the presence of
a.c. spectrum.

Consider a general solution
\[
x_{n+1}+v_nx_n+x_{n-1}=zx_n,\quad n=1,2,\ldots
\]
\[
X_n=\left[
\begin{array}{c}
x_{n+1}\\ x_n
\end{array}
\right], \quad n=0,1,\ldots
\]

\[
X_{n}=(\Omega+V_n) X_{n-1}, \quad \Omega=\left[
\begin{array}{cc}
z & -1\\
1 & 0
\end{array}
\right],\quad V_n=v_n\left[
\begin{array}{cc}
-1 & 0\\
0 & 0
\end{array}
\right],\quad n=1,2,\ldots
\]
If $Z_m=X_{mq}$, then
\begin{equation}
Z_{m+1}=T_mZ_m, \quad m=0,1,\ldots \label{recz}
\end{equation}
where
\[
T_m=(\Omega+V_{mq+q})\ldots (\Omega+V_{mq+1})
\]

Let
\[
k^{\pm 1}(z)=\frac{z\mp\sqrt{z^2-4}}{2}
\]
so that $k(z)$ maps $\mathbb{C}\backslash [-2,2]$ onto $\mathbb{D}$
conformally.

Notice that if $\tilde{P},P$ and $\tilde{Q},Q$ are $(q, q-1)$-th
and the $(q-1, q-2)$-th polynomials corresponding to the first
Jacobi coefficients $v_1, \ldots, v_q$ and $(1,0)^t$, $(0,1)^t$
initial conditions, then
\[
T_m=\left[
\begin{array}{cc}
\tilde{P}_m & \tilde{Q}_m\\
P_m & Q_m
\end{array}
\right]
\]
Notice that $\det T_m=1$ and therefore
$\tilde{P}_mQ_m-P_m\tilde{Q}_m=1$. We also have
\[
\Omega^q=\left[
\begin{array}{cc}
\displaystyle \frac{k^{q+1}-k^{-q-1}}{k-k^{-1}} &
\displaystyle-\frac{k^{q}-k^{-q}}{k-k^{-1}}\\\\
\displaystyle \frac{k^{q}-k^{-q}}{k-k^{-1}} & \displaystyle
-\frac{k^{q-1}-k^{-q+1}}{k-k^{-1}}
\end{array}
\right]
\]
and therefore
\begin{equation}
\tilde{P}_m+Q_m=k^q+k^{-q}+d(k,v_{mq+1}, \ldots,
v_{mq+q}),\label{petshop}
\end{equation}
The function $d(\cdot)$ is a polynomial in $v_{mq+1}, \ldots,
v_{mq+q}$ and
\begin{equation}
d(k,v_{mq+1}, \ldots, v_{mq+q})\to 0, \,
m\to\infty\label{petshop1}
\end{equation}

Introduce $\lambda_{1(2)}^{(m)}$ by
\begin{equation}
\lambda_{1(2)}^{(m)}=\frac{\tilde{P}_m+Q_m\mp
\sqrt{(\tilde{P}_m+Q_m)^2-4}}{2} \label{roots}
\end{equation}
These are the eigenvalues of $T_m$.

%\[
%=\frac{k^q+k^{-q}\mp
%\sqrt{\displaystyle(k^q-k^{-q})^2+d_m(k,v)}}{2},\,\, d(k,v_{mq+1},
%\ldots, v_{mq+q})\to 0,{\rm \,\,as}\,\, m\to \infty
%\]
Let us take $U_m$
\begin{equation}
U_m=\left[
\begin{array}{cc}
-Q_m+\lambda_1^{(m)} & -Q_m+\lambda_2^{(m)}\\
P_m & P_m
\end{array}
\right]\label{um}
\end{equation}
Then we have
\begin{equation}
U_{m+1}^{-1}U_m=\frac{1}{P_{m+1}(\lambda_1^{(m+1)}-\lambda_2^{(m+1)})}\label{ratu}
\end{equation}
\[
\times\left[
\begin{array}{cc}
P_{m+1}(\lambda_1^{(m)}-Q_m)-P_m(\lambda_2^{(m+1)}-Q_{m+1}), & P_{m+1}(\lambda_2^{(m)}-Q_m)-P_m(\lambda_2^{(m+1)}-Q_{m+1})\\
-P_{m+1}(\lambda_1^{(m)}-Q_m)+P_m(\lambda_1^{(m+1)}-Q_{m+1}),&-P_{m+1}(\lambda_2^{(m)}-Q_m)+P_m(\lambda_1^{(m+1)}-Q_{m+1})
\end{array}\right]
\]
The matrix $U_m$ can be used to diagonalize $T_m$ as follows
\[
T_m(z)=U_m(z)\Lambda_m(z)U_m^{-1}(z)
\]
where
\[
\Lambda_m(z)=\left[\begin{array}{cc}
\lambda_1^{(m)}(z) & 0\\
0& \lambda_2^{(m)}(z)
\end{array}\right]
\]
\begin{lemma}
The matrix $T_m$ has eigenvalues $\lambda_{1(2)}^{(m)}(z)$ such
that
\begin{equation}
\lambda_{1}^{(m)}(z)\cdot \lambda_{2}^{(m)}(z)=1, \quad z\in
\mathbb{C} \label{unif}
\end{equation}
\[
|\lambda_{1(2)}^{(m)}(z)|=1, \quad z\in
[z_j+\delta(v),z_{j+1}-\delta(v)], \quad  j=0,\ldots, q-1
\]
and
\[
\delta(v)\to 0
\]
as
\[\zeta_m=\max\{|v_{mq+q}|, \ldots ,|v_{mq+1}|\}\to 0
\] Here $z_j=2\cos(\pi -\pi j/q), \quad j=0, \ldots ,q$.\label{lemmal1}
\end{lemma}
\begin{proof}
The first identity follows from $\det T_m(z)=1$. Since the function
$d(\cdot)$ in (\ref{petshop}) is real for real $z$, the second one
is immediate from (\ref{petshop}), (\ref{petshop1}), and
(\ref{roots}).
\end{proof}
Notice that $\Omega^q$ has eigenvalues
\[
\omega_{1(2)}=k^{\mp q}(z)=\left(\frac{z\pm
\sqrt{z^2-4}}{2}\right)^q
\]
We  have the following elementary perturbation result
\begin{lemma}
If \, $0\leq \Im z\leq  1, \, z_j+\delta<\Re z<z_{j+1}-\delta, \,
j=0,\ldots, q-1 $, then
\begin{eqnarray}
\begin{array}{c}
\lambda_{1(2)}^{(m)}(z)=\omega_{1(2)}(z)+\underline{O}(\zeta_m)\\\\
\sum\limits_{m=0}^\infty
\left|\lambda_{1(2)}^{(m+1)}(z)-\lambda_{1(2)}^{(m)}(z)\right|^2
<C(\delta)
\end{array}\label{dif}
\end{eqnarray}
and $\delta>0$ is fixed arbitrarily small number. Moreover, for all
$z$ in these domains we have the following estimate
\begin{equation}
 \ln |\lambda_1^{(m)}(z)|=(C+\underline{O}(\zeta_m))\Im z
 \label{unif-exp}
\end{equation}
with some positive constant $C$. We also assume here that
$m>m_0(\delta)$ and $m_0(\delta)$ is large depending on $\delta$.
 \label{ratiol}
\end{lemma}
\begin{proof}
For $d(\cdot)$, we have $|d(k,v_{mq+1}, \ldots,
v_{mq+q})|<C|\zeta_m|$. It is also a polynomial in $v_{mq+1},
\ldots, v_{mq+q}$. Then, (\ref{dif}) follows from the Mean Value
Theorem and (\ref{petshop}), (\ref{petshop1}), and (\ref{roots}). To
prove (\ref{unif-exp}), we fix $j$ and consider the following
function
\[
h(z)=\ln |\lambda_1^{(m)}(z)/\omega_1(z)|
\]
harmonic in the domain of interest:\quad $0\leq \Im z\leq  1, \,
z_j+\delta<\Re z<z_{j+1}-\delta$. On the real line, i.e., for\,
$z_j+\delta<\Re z<z_{j+1}-\delta,\, \Im z=0$, we have $h(z)=0$ and
at all other points we have
\[
h(z)=\underline{O}(\zeta_m)
\]
Therefore, the interpolation lemma \ref{interpol} gives
\[
h(z)=\underline{O}(\zeta_m)\Im z
\]
For $\omega_1(z)$ we have $|\omega_1(z)|>1+C\Im z$ with positive
$C$, and one gets (\ref{unif-exp}).
\end{proof}
Let us find $Z_m$ in the form $Z_m=U_m S_m$ and
\begin{equation}
S_{m+1}=U_{m+1}^{-1}T_mU_mS_m= U_{m+1}^{-1} U_m \Lambda_m S_m=
[U_{m+1}^{-1}(U_m-U_{m+1})+I]\Lambda_m S_m \label{s-rec}
\end{equation}

\begin{lemma}
If \, $0\leq \Im z\leq  1, \, z_j+\delta<\Re z<z_{j+1}-\delta, \,
j=0,\ldots, q-1 $, then for the matrix norms we have
\[
\|U_{m+1}^{-1}(U_m-U_{m+1})\|_{\ell^2}<C
\]\label{difl}
\end{lemma}
\begin{proof}
Away from the points $z_j$, $\|U_{m+1}^{-1}\|$ is uniformly bounded
and the proof follows immediately from (\ref{dif}) and (\ref{um}).
\end{proof}

We need the following
\begin{theorem}
Let \[\Psi_{n+1}=(I+W_n)\left[
\begin{array}{cc}
\kappa_n & 0\\
0 & \kappa_n^{-1}
\end{array}
\right]\Psi_n, \,W_n=\left[
\begin{array}{cc}
\alpha_n & \beta_n\\
\gamma_n & \delta_n
\end{array}
\right], \, \Psi_0= \left[
\begin{array}{c}
1\\
0
\end{array}
\right], \] where $\kappa_n\in \mathbb{C}$,
$C>|\kappa_n|>|\kappa|>1$, the sequence $\zeta_n=\|W_n\|\in
\ell^2(\mathbb{Z}^+)$ and its $\ell^2$ norm is small. Assume also
that there is a constant $0\leq \upsilon<1$ so that
\begin{equation}
  \left|\,\ln \prod\limits_{n=k}^l \left|1+\alpha_n\right|\right|\leq C+\upsilon \sqrt{l-k},\quad
\left|\,\ln \prod\limits_{n=k}^l
\left|1+\delta_n\right|\right|\leq C+\upsilon\sqrt{l-k}\label{ups}
\end{equation}
Then,
\[
\Psi_n= p_n\left[
\begin{array}{c}\phi_n
\\
\nu_n
\end{array}
\right],
\]
where
\[
p_n=\prod\limits_{j=0}^{n-1} \kappa_j(1+\alpha_j), \quad p_0=1
\]
and
\begin{equation}
|\phi_n|, |\nu_n|\leq
C\exp\left(\frac{C}{|\kappa|-1}\exp\left[\frac{C\upsilon^2}{|\kappa|-1}\right]\right)\label{one}
\end{equation}
Moreover, for any fixed $\epsilon>0$ and any $\kappa:
|\kappa|>1+\epsilon$, we have
\begin{equation}
|\phi_n|>C>0, \quad |\nu_n|<C\|\zeta\|_2 \label{pepsi}
\end{equation}
uniformly in $n$ provided that $\|\zeta\|_2$ is small enough.
\label{theorem21}
\end{theorem}
\begin{proof}
Let
\[
\Psi_n=\left[
\begin{array}{cc}
p_n & 0\\
0 & q_n
\end{array}
\right]Y_n
\]
where
\[
q_n= \prod\limits_{j=0}^{n-1} \kappa_j^{-1}(1+\delta_j), \, q_0=1
\]
Then,
\[
Y_{n+1}=\left[
\begin{array}{cc}
1 & q_np_{n+1}^{-1}\kappa_n^{-1}\beta_n\\
p_nq_{n+1}^{-1}\kappa_n\gamma_n & 1
\end{array}
\right]Y_n
\]
If
\[
Y_n=\left[\begin{array}{c}
\phi_n\\
v_n
\end{array}
\right]
\]
then
\begin{equation}
\nu_n=q_np_n^{-1}v_n \label{rel17}
\end{equation}
and we have the following equations
\[
\phi_n=1+\sum\limits_{j=0}^{n-1} q_jp_{j+1}^{-1}\kappa_j^{-1}\beta_j
v_j
\]
\begin{equation}
v_n=\sum\limits_{j=0}^{n-1} p_jq_{j+1}^{-1}\kappa_j \gamma_j\phi_j
\label{vvv}
\end{equation}
For $\phi_n$:
\[
\phi_n=1+\sum\limits_{k=0}^{n-2} \phi_k \epsilon_{k,n},\quad
\epsilon_{k,n}=\gamma_k p_kq_{k+1}^{-1}\kappa_k
\sum\limits_{j=k+1}^{n-1}\kappa_j^{-1}q_jp_{j+1}^{-1}\beta_j
\]
For $\epsilon_{k,n}$,
\begin{equation}
|\epsilon_{k,n}|\leq \epsilon_k=C |\gamma_k|
\sum\limits_{j=k+1}^{\infty}|\beta_j| \cdot |\kappa|^{-2(j-k)}\cdot
\prod\limits_{l=k}^{j} \left|\frac{1+\delta_l}{1+\alpha_l}\right|
\label{seqe}
\end{equation}
\[
<C |\gamma_k| \sum\limits_{j=k+1}^{\infty}|\beta_j|
|\kappa|^{-2(j-k)} \cdot \exp\left( C\upsilon \sqrt{j-k}\right)
\]
From the obvious inequality
\[
\sum_{j=0}^\infty |\kappa|^{-j} \cdot \exp\left( C\upsilon
\sqrt{j}\right)<\frac{C}{\ln
|\kappa|}\exp\left[\frac{C\upsilon^2}{\ln |\kappa|}\right]
\]
and Young's inequality for convolutions, we get an estimate for the
$\ell^1$ norm of the sequence $\epsilon$ introduced in (\ref{seqe}):
\begin{equation}
\|\epsilon\|_1\leq \frac{C}{|\kappa|-1}\cdot
\exp\left[\frac{C\upsilon^2}{|\kappa|-1}\right]\cdot
\|\gamma\|_2\cdot \|\beta\|_2
\end{equation}
This yields the same estimates for
\[
\sum\limits_{k=0}^n |\epsilon_{k,n}|
\]
uniformly in $n$. Now, to prove (\ref{one}), one can use the
following lemma below.
\begin{lemma}
If $\, x_n, v_n\geq 0, \, x_0=1$, and
\[
x_{n+1}\leq\sum\limits_{j=0}^n v_j x_j
\]
for all $n>0$, then
\begin{equation}
x_n\leq v_0\exp\left[\sum\limits_{j=1}^{n-1} v_j\right],\,\, n\geq
2;\quad\quad x_1\leq v_0 \label{gb}
\end{equation}
\end{lemma}
\begin{proof}\footnote{This lemma can also be proved by
induction.} Consider the functions
\[
x(t)=\sum\limits_{j=0}^\infty x_j \chi_{[j,j+1)}(t),\quad
v(t)=\sum\limits_{j=0}^\infty v_j \chi_{[j,j+1)}(t)
\]
We have
\[
x(t)\leq v_0+\int\limits_1^{[t]-1} x(s)v(s)ds\leq
v_0+\int\limits_1^t x(s)v(s)ds,\quad t>1
\]
The application of Gronwall-Bellman inequality gives (\ref{gb}).
\end{proof}
The estimate for $\nu_n$ and the line (\ref{pepsi}) are
straightforward corollaries from the bound for $\|\phi\|_\infty$,
(\ref{rel17}), (\ref{vvv}), and the Cauchy-Schwarz inequality.
\end{proof}

Introduce
\begin{equation}
U_{m+1}^{-1}(U_m-U_{m+1})=\left[
\begin{array}{cc}
\alpha_m & \beta_m\\
\gamma_m & \delta_m
\end{array}
\right]\label{not}
\end{equation}
Then, (\ref{s-rec}) can be rewritten as
\begin{equation}
S_{m+1}=\left[I+\left[
\begin{array}{cc}
\alpha_m & \beta_m \\
\gamma_m & \delta_m
\end{array}
\right]\right]\times\left[
\begin{array}{cc}
\lambda_1^{(m)} & 0 \\
0 & \lambda_2^{(m)}
\end{array}
\right]    S_m \label{lrec1}
\end{equation}

Now, let us apply the theorem \ref{theorem21}  to (\ref{lrec1}). For
each $k=0,\ldots, q-1$, consider $z$ in the following domain: $0\leq
\Im z<1$, $z_k+\delta<\Re z<z_{k+1}-\delta$ where $\delta$ is a
small positive number. We have $\lambda_1^{(m)}\cdot
\lambda_2^{(m)}=1$ and $|\lambda_1^{(m)}|>(1+C\Im z)$ by lemma
\ref{lemmal1} and (\ref{unif-exp}).
 In our
notations
\[
W_m=\left[
\begin{array}{cc}
\alpha_m & \beta_m \\
\gamma_m & \delta_m
\end{array}
\right]
\]
and $\kappa_m=\lambda_1^{(m)}$. The estimate $\|W_m\|\in
\ell^2(\mathbb{Z}^+)$ follows from the lemma \ref{difl}. Now, let us
control
\[
\prod\limits_{j=0}^n (1+\alpha_n), \quad \prod\limits_{j=0}^n
(1+\delta_n)
\]
and the constant $\upsilon(z)$ in (\ref{ups}).

\begin{theorem} For
$z: 0\leq \Im z<1,\, z_k+\delta<\Re z<z_{k+1}-\delta,\,
k=0,1,\ldots, q-1$, we have
\begin{equation}
  \left|\,\ln \prod\limits_{n=k}^l \left|1+\alpha_n\right|\right|\leq C+(C\Im z)  (l-k)^{1/2},\quad
\left|\,\ln \prod\limits_{n=k}^l
\left|1+\delta_n\right|\right|\leq C+(C\Im z)(l-k)^{1/2}
\end{equation}
as long as $\|v\|_\infty$ is small.\label{theorem22}
\end{theorem}
\begin{proof}
Consider the product for $\alpha_n$, the product for $\delta_n$
can be treated in the same way.  In notations of (\ref{ratu}),
(\ref{not}),
\[
1+\alpha_m(z)=\frac{P_{m+1}
(\lambda_1^{(m)}-Q_m)-P_m(\lambda_2^{(m+1)}-Q_{m+1})}{P_{m+1}(\lambda_1^{(m+1)}-\lambda_2^{(m+1)})}=
1+t_m^1+t_m^2+t_m^3+t_m^4+t_m^5
\]
where
\[
t_m^1=\frac{P_{m}Q_{m+1}-P_{m+1}Q_m}
{P_{m+1}(\lambda_1^{(m+1)}-\lambda_2^{(m+1)})}
\]
\[
t_m^2=-\frac{P_{m+1}-P_m}{2P_{m+1}}
\]
\[
t_m^3=\frac{(P_{m+1}-P_m)(\tilde{P}_{m+1}+Q_{m+1})}{2P_{m+1}(\lambda_1^{(m+1)}-\lambda_2^{(m+1)})}
\]
\[
t_m^4=-\frac{\tilde{P}_{m+1}-\tilde{P}_m+Q_{m+1}-Q_m}{2(\lambda_1^{(m+1)}-\lambda_2^{(m+1)})}
\]
\[
t_m^5=-\frac{(\lambda_1^{(m+1)}-\lambda_2^{(m+1)})-(\lambda_1^{(m)}-\lambda_2^{(m)})}
{2(\lambda_1^{(m+1)}-\lambda_2^{(m+1)})}
\]
since
\[
\frac{P_{m+1}\lambda_1^{(m)}-P_m\lambda_2^{(m+1)}}{{P_{m+1}(\lambda_1^{(m+1)}-\lambda_2^{(m+1)})}}
=\frac{P_{m+1}(
\lambda_1^{(m)}-\lambda_2^{(m+1)})}{{P_{m+1}(\lambda_1^{(m+1)}-\lambda_2^{(m+1)})}}+
\]
\[
+\frac{(P_{m+1}
-P_m)\lambda_2^{(m+1)}}{{P_{m+1}(\lambda_1^{(m+1)}-\lambda_2^{(m+1)})}}
\]
\[
\frac{
\lambda_1^{(m)}-\lambda_2^{(m+1)}}{{\lambda_1^{(m+1)}-\lambda_2^{(m+1)}}}=1-\frac{
\lambda_1^{(m+1)}-\lambda_1^{(m)}}{{\lambda_1^{(m+1)}-\lambda_2^{(m+1)}}}
\]
\[
\frac{
\lambda_1^{(m+1)}-\lambda_1^{(m)}}{{\lambda_1^{(m+1)}-\lambda_2^{(m+1)}}}=\frac{(\lambda_1^{(m+1)}-\lambda_2^{(m+1)})-(\lambda_1^{(m)}-\lambda_2^{(m)})}
{2(\lambda_1^{(m+1)}-\lambda_2^{(m+1)})}+
\]
\[
+\frac{\tilde{P}_{m+1}-\tilde{P}_m+Q_{m+1}-Q_m}{2(\lambda_1^{(m+1)}-\lambda_2^{(m+1)})}
{\rm \quad as\,\,  follows\,\,  from\,\, (\ref{roots})}
\]
and
\[
\frac{(P_{m+1}
-P_m)\lambda_2^{(m+1)}}{{P_{m+1}(\lambda_1^{(m+1)}-\lambda_2^{(m+1)})}}=-\frac{P_{m+1}-P_m}{2P_{m+1}}
+\frac{(P_{m+1}-P_m)(\tilde{P}_{m+1}+Q_{m+1})}{2P_{m+1}(\lambda_1^{(m+1)}-\lambda_2^{(m+1)})}
\]

 Obviously, all $t_m^j\in \ell^2(\mathbb{Z}^+), j=1,\ldots, 5$ and are small
if $\|v\|_\infty$ is small.

Therefore, we have
\begin{equation}
  \left|\,\Re \ln \prod\limits_{n=k}^l \left(1+\alpha_n\right)\right|\leq
  \left|\sum_{j=1}^5 \Re \sum\limits_{n=k}^l t_n^j\right| +C
\end{equation}
We need the following lemma
\begin{lemma}
If $\epsilon_n\to \epsilon\in \mathbb{C}$,
$|\epsilon_{n+1}-\epsilon_n|\in \ell^2$, and $f(z)$ is holomorphic
around $\epsilon$. Then,
\[
\left|\sum\limits_{n=k}^l
(\epsilon_{n+1}-\epsilon_n)f(\epsilon_n)\right|<C
\]
for $k$ and $l$ large enough.\label{lbp}
\end{lemma}
\begin{proof}
Consider $g(z)$ holomorphic around $\epsilon$ such that
$g(\epsilon)=0$ and $g'(z)=f(z)$. The Taylor formula around
$\epsilon_n$ yields
\[
g(\epsilon_{l+1})-g(\epsilon_k)=\sum\limits_{n=k}^l
g(\epsilon_{n+1})-g(\epsilon_n) =\sum\limits_{n=k}^l
\left[(\epsilon_{n+1}-\epsilon_n)f(\epsilon_n)+\underline{O}(|\epsilon_{n+1}-\epsilon_n|^2)\right]
\]
which shows that the second term in the right hand side is bounded.
\end{proof}
Taking $f(z)=z^{-1}$ in the lemma yields
\[
\left|\sum\limits_{n=k}^l t_n^j\right|<C
\]
where $j=2,5$. Now, for the other $j$, consider the following
functions
\[\Re \sum\limits_{n=k}^l t_n^j(z)
\]
which are harmonic near the intervals of interest. By
Cauchy-Schwarz, we have
\begin{equation}
\left|\Re \sum\limits_{n=k}^l
t_n^j(z)\right|<C(l-k)^{1/2}\label{prop1}
\end{equation}
for any $z$ in the specified domain. For real $z$ within the
intervals $[z_k+\delta,z_{k+1}-\delta]$,
\begin{equation}
\Re \sum\limits_{n=k}^l t_n^j(z)=0 \label{prop2}
\end{equation}
Indeed, all polynomials $P,\tilde{P},Q,\tilde{Q}$ are real for real
$z$. On the other hand, for real $z$ within these intervals and
small $\|v\|_\infty$, we have
$\lambda_1^{(m)}=\overline{\lambda_2^{(m)}}$, and so
$\lambda_1^{(m)}-\lambda_2^{(m)}$ is purely imaginary.  Now, the
theorem follows from (\ref{prop1}), (\ref{prop2}) and interpolation
lemma \ref{interpol}.

\end{proof}

As a corollary, we get
\begin{theorem}
Consider $z: 0\leq \Im z<1,\, z_k+\delta<\Re z<z_{k+1}-\delta,\,
k=0,1,\ldots, q-1$ and introduce $\kappa_j(z)$ and $\alpha_j(z)$ as
before. If $v_1=v_2=\ldots=v_q=0$ and $X_0=(k^{-1}(z),1)^t$, then we
have the following estimates for the solution of the Schr\"odinger
recursion:
\[
|x_{nq}(z)|\leq \left|\prod\limits_{j=0}^{n-1}
{\kappa_j(z)(1+\alpha_j(z))} \right|\cdot \exp\left(\frac{C}{\Im
z}\right)
\]
where the first factor in the r.h.s. is uniformly bounded for $z\in
(z_k+\delta,z_{k+1}-\delta)$, $k=0,1,\ldots,q-1$. Moreover, for any
fixed $\epsilon>0$ and any $z: \Im z>\epsilon>0$, we have
\[
\left|x_{nq}(z) \left({\prod\limits_{j=0}^{n-1}
\kappa_j(z)(1+\alpha_j(z))}\right)^{-1}\right|>C>0
\]
uniformly in $n$, provided that $\|v\|_\infty+\|\delta^{(q)}v\|_2$
is small enough. \label{est-key}
\end{theorem}
\begin{proof}
We have
\[
U_0^{-1}=\frac{1}{k^{-q}-k^q} \left[
\begin{array}{cc}
1& -k\\
-1 & k^{-1}
\end{array}
\right]
\]
If
\begin{equation}
X_{mq}=Z_m=U_mS_m \label{ts}\end{equation} then we have recursion
(\ref{s-rec}) for $S_m$ and $
S_0=(k(z)-k^{-1}(z))(k^q(z)-k^{-q}(z))^{-1}(1,0)^t $. From
(\ref{ts}), an explicit form (\ref{um}) for $U_m$, and theorems
\ref{theorem21}, \ref{theorem22}, we get the statement of the
theorem.
\end{proof}
{\bf Remark.} One can easily see that the analysis suggested above
proves the asymptotics for $x_{n}(z)$ in the corresponding domains
of complex plane, not just bounds from above and below. It is also
conceivable that the complicated product of $1+\alpha_j$ in the
asymptotics can be simplified and possibly eliminated (as for $q=1$)
by applying some analog of lemma~\ref{lbp}. We do not pursue it here
and employ technique which is somewhat more powerful.

Now, we are ready to prove the main result of this paper.
\begin{theorem}
If the Jacobi matrix $J$ given by (\ref{jacobi}) has coefficients
$v_n\to 0$ and $v_{n+q}-v_n\in \ell^2(\mathbb{Z}^+)$ for some $q$,
then $\sigma_{ac}(J)=[-2,2]$.
\end{theorem}
\begin{proof}
By Weyl's theorem, the essential spectrum is $[-2,2]$. Take any
small $\delta>0$. We will show that all  intervals $[z_j+\delta,
z_{j+1}-\delta], \quad j=0,\ldots, q-1$ are  contained in the
support the a.c. spectrum of $J$. Fix $\delta>0$. Just like in the
theorem \ref{th1}, we can assume that $\|v\|_\infty+\|\delta^{(q)}
v\|_2$ is as small as we wish.

Then, we will prove

\begin{equation} \int\limits_{z_j+\delta}^{z_{j+1}-\delta} \ln
\rho'(z)dz>-\infty, \quad j=0,\ldots, q-1
\label{inter}\end{equation} where $\rho(z)$ is the spectral measure
of the matrix $J$. Consider the truncated potential $v^N=v\cdot
\chi_{j<N}$, the corresponding matrix $J^N$, and the spectral
measure $\rho^N$, (we take $N=qm$). Since the potential is finitely
supported, there exists the Jost solution:
\[
x_{n+1}^N+v^N_nx_n^N +x_{n-1}^N=zx_n^N
\]
with the following asymptotics
\[
x_n^N=k^n, \quad n>N
\]
Then, the  factorization (see, \cite{ks}, (1.32))
\[
\pi\cdot(\rho^N)'(z)=\frac{\sqrt{4-z^2}}{|x_0^N(z)|^2},\quad z\in
[-2,2]
\]
holds.

Obviously, $x_n^N$ is the solution of the Cauchy problem
$x_N^N=k^N$, $x_{N+1}^N=k^{N+1}$ and $x_0^N$ can be obtained by
solving the recursion ``backward". The theorem~\ref{est-key} can be
applied then. Introduce the function (modified Jost function)
\[
f_N(z)=x_0^N(z)\left(k^N\prod\limits_{j=0}^{m-1}
(1+\alpha_j(z))\kappa_j(z)\right)^{-1}
\]
where $\alpha_j$ and $\kappa_j$ are taken from the theorem
\ref{est-key} (with ``backward" ordering for potential). From this
theorem,
\[
|f_N(z)|<\exp\left(\frac{C}{\Im z}\right)
\]
uniformly in $N$ as long as $z_j+\delta<\Re z<z_{j+1}-\delta, 0<\Im
z<1$. On the real line,

\[
\left|k^N\prod\limits_{j=0}^{m-1}
(1+\alpha_j(z))\kappa_j(z)\right|\sim 1
\]
and we obtain
\begin{equation}
\int\limits_{z_j+\delta}^{z_{j+1}-\delta}\ln
(\rho^N)'(z)dz>-C_1-C_2\int\limits_{z_j+\delta}^{z_{j+1}-\delta}\ln
|f_N(z)|dz \label{aa1}
\end{equation}
\begin{equation}
\int\limits_{z_j+\delta}^{z_{j+1}-\delta}\left(-\ln
|f_N(z)|\right)^{+}dz<C_1\int\limits_{z_j+\delta}^{z_{j+1}-\delta}\left(\ln
(\rho^N)'(z)\right)^+dz+C_2<C\label{aa2}
\end{equation}
uniformly in $N$. Moreover,
\[
|f_N(z)|>C>0
\]
if $\Im z>\epsilon(v)$, also uniformly in $N$.

 The function $-\ln |f_N(z)|$ is harmonic in
$z_j+\delta<\Re z<z_{j+1}-\delta, 0<\Im z<1$ because $x_0(z)$ has no
zeroes in $\mathbb{C}^+$ (otherwise the asymptotics at infinity
would yield the complex eigenvalue for $J^N$ which is impossible)
and $1+\alpha_j(z)\neq 0$ since $\|v\|_\infty$ is small. This
function is also continuous up to the boundary since the potential
is finitely supported. Due to (\ref{aa2}), the lemma \ref{lemma2} of
Appendix is then applicable. It yields
\[\int\limits_{z_j+\delta}^{z_{j+1}-\delta}(-\ln|f_N(z)|)^{-} >-C\quad j=0,\ldots, q-1
\]
uniformly in $N$. By (\ref{aa1}), we also have
\[\int\limits_{z_j+\delta}^{z_{j+1}-\delta}\ln
(\rho^N)'(z)dz>C, \quad j=0,\ldots, q-1
\]
Now, notice that $d\rho^{N}\to d\rho$ weakly as $N\to\infty$ (see,
e.g., \cite{DK}) and the semicontinuity of the entropy argument from
\cite{ks} (see corollary 5.3) gives (\ref{inter}). Since $\delta>0$
was arbitrary positive, that proves that the a.c. spectrum is
supported on $[-2,2]$.
\end{proof}

\bigskip

 {\bf Acknowledgements.} The author is grateful to A. Kiselev and B. Simon
 for useful discussion. He also wants to thank Uri Kaluzhny and Mira
 Shamis for pointing out one mistake in the original version of
 this manuscript.
 This work was supported by Alfred P.
Sloan Research Fellowship and NSF Grant DMS-0500177.

\section{Appendix}

In this Appendix, we prove several auxiliary statements used in the
main text.
\begin{lemma}
Assume that $f(z)$ is harmonic in $\Pi=a<\Re z<b, \, 0<\Im z<c$ and
is continuous on the closure $\overline{\Pi}$. Then, two estimates
\[
|f(z)|<\gamma, \quad z\in \overline\Pi
\]
and \[f(z)=0,\quad z\in [a,b]
\]
imply
\[
|f(z)|\leq C\gamma\Im z,\quad a+\delta<\Re z<b-\delta
\]
and the constant $C$ depends on the domain $\Pi$ and $\delta$
only.\label{interpol}
\end{lemma}
\begin{proof}
Map $\Pi$ conformally onto the upper half-plane such that $[a,b]$
goes to, say, $[-1,1]$. Let $g(\xi)$ be the transplantation of
$f(z)$. Then, an application of the Poisson representation
\[
g(x+iy)=\frac{y}{\pi}\int_{\mathbb{R}} \frac{g(t)}{(t-x)^2+y^2}dt
\]
provides the necessary estimate.
\end{proof}

 Consider the domain $\Sigma=\{z=x+iy: a<x<b, 0<y<c\}$ in
the complex plane.  The next lemma is rather simple and is taken
from \cite{killip}. Below, the constants $C$ are all positive,
$f$--independent, and can change from one formula to another;
$f^{\pm}$ denotes the positive/negative parts of a real--valued
function $f$ (i.e., $f^+=f$ if $f\geq 0$ and $f=0$ otherwise,
$f^-=f-f^+$).

\begin{lemma}
Assume $f(z)$ is harmonic on $\Sigma$, continuous on
$\overline{\Sigma}$, and
\begin{equation}
\int\limits_{a}^{b} f^+(x)dx<C \label{cond1}
\end{equation}
\begin{equation}
f(z)<Cy^{-\alpha}, \quad z\in \Sigma,\quad \alpha>0 \label{cond2}
\end{equation}
and
\begin{equation}
f(z_1)>-C \label{cond3}
\end{equation}
for some $z_1\in\Sigma$. Then,
\[
\int\limits_{a+\delta}^{b-\delta} f^-(x)dx>-C(\delta)
\]

\end{lemma}
\begin{proof}
Take isosceles triangle $ABC$ as shown on Figure 1 below and write
the mean value formula for $f$
\[
\int\limits_{AB} f(z)\omega(z)d|z|+\int\limits_{AC}
f(z)\omega(z)d|z|+\int\limits_{BC} f(z)\omega(z)d|z|=f(z_1)
\]
where $\omega$ is the harmonic measure for triangle and $z_1$
assuming that $z_1$ is inside it.  That can be rewritten as
\[
\int\limits_{AB} f^-(z)\omega(z)d|z|+\int\limits_{AC}
f^-(z)\omega(z)d|z|+\int\limits_{BC} f^-(z)\omega(z)d|z|=
\] \[=f(z_1)-\int\limits_{AB}
f^+(z)\omega(z)d|z|-\int\limits_{AC}
f^+(z)\omega(z)d|z|-\int\limits_{BC} f^+(z)\omega(z)d|z|
\]
By taking the angle $\gamma$ small enough we can always guarantee
that $\omega$ will decay fast at $A$ and $B$ so that the integrals
of $f^+$ over $AC$, $BC$ are uniformly bounded due to
(\ref{cond2}). The integral of $f^+$ over $AB$ is  bounded due to
(\ref{cond1}). Therefore, since $f(z_1)$ is bounded from below by
(\ref{cond3}), we get the uniform estimate for the left-hand side.
If $z_1$ is not inside the triangle, an analogous argument would
work if one takes a slightly different shape (e.g., a curved
triangle).
\end{proof}
In the main text, the condition that $z_1$ is inside a triangle
can always be realized by assuming the norm of  $v$ to be small.

\unitlength=0.3mm
\begin{picture}(300,200)
 \put(50,100){\line(1,0){300}}
\put(50,100){\line(4,1){150}}\put(202, 139){$C$}
\put(350,100){\line(-4,1){150}}
 \put(200,120){\circle*{2}}
\put(205,118){$z_1$} \qbezier[1200](100,100)(99,106)(95,111)
\put(102,103){$\gamma$} \put(50,90){$a$} \put(45, 102){$A$}
\put(350,90){$b$} \put(350, 102){$B$}

 \put(180,50){Figure 1}
\end{picture}

The slight modification of the argument provides the next lemma
which is crucial for our considerations. It allows to obtain the
entropy bounds from very rough estimates on modified Jost
functions.
\begin{lemma}
Assume $f(z)$ is harmonic on $\Sigma$, continuous on
$\overline{\Sigma}$, and
\[
\int\limits_{a}^{b} f^+(x)dx<C
\]
\begin{equation}
f(z)>-Cy^{-\alpha},\quad z\in \Sigma, \quad\alpha>0 \label{l2est2}
\end{equation}
and
\begin{equation}
f(z)<C \label{l2est3}
\end{equation}
for $\Im z>d(\alpha)=C(1+\alpha)^{-1}>0$. Then,
\[
\int\limits_{a+\delta}^{b-\delta} f^-(x)dx>-C(\delta)
\]\label{lemma2}
\end{lemma}
\begin{proof}
Take two equal triangles $ABC$ and $ADC$ as shown on the Figure 2.
As in the previous lemma,
\[
\int\limits_{AB} f^-(z)\omega(z)d|z|+\int\limits_{AC}
f^-(z)\omega(z)d|z|+\int\limits_{CB} f^-(z)\omega(z)d|z|=
\]
\begin{equation}
=f(z_1)-\int\limits_{AB} f^+(z)\omega(z)d|z|-\int\limits_{AC}
f^+(z)\omega(z)d|z|-\int\limits_{CB} f^+(z)\omega(z)d|z|
\label{interm}
\end{equation}
By the same arguments, we only need to provide a uniform estimate
from above for the integrals over $AC$ and $BC$. Consider the
integral over $AC$, the other side can be treated in the same way.
Consider $ACD$ with $z_2$ being symmetric to $z_1$ with respect to
the line $AC$. The harmonic measure for $ACD$ and $z_2$ is equal
to the same $\omega$ and we have

\[
\int\limits_{AD} f^+(z)\omega(z)d|z|+\int\limits_{AC}
f^+(z)\omega(z)d|z|+\int\limits_{DC} f^+(z)\omega(z)d|z|=
\] \[=f(z_2)-\int\limits_{AD}
f^-(z)\omega(z)d|z|-\int\limits_{AC}
f^-(z)\omega(z)d|z|-\int\limits_{DC} f^-(z)\omega(z)d|z|
\]
The estimate (\ref{l2est3}) gives a bound for $f(z_2)$ from above
and $(\ref{l2est2})$ controls the other terms in the right-hand
side if the angle $\gamma$ is small enough. Thus,
\[
\int\limits_{AC} f^+(z)\omega(z)d|z|<C
\]
and that finishes the proof due to (\ref{interm}). Choosing
$d(\alpha)=C(1+\alpha)^{-1}$ with small $C$ guarantees that
(\ref{l2est3}) is applicable to $z_2$.
\end{proof}

\begin{picture}(300,300)
\put(45, 102){$A$}\put(350, 102){$B$}
 \put(50,100){\line(1,0){300}} \put(50,100){\line(4,1){150}}
 \put(209, 137){$C$} \put(318, 238){$D$}
\put(350,100){\line(-4,1){150}}
\qbezier[1200](50,100)(315,241)(315,241)
\qbezier(200,138)(200,138)(315,241) \put(200,120){\circle*{2}}
\put(205,118){$z_1$} \put(190,153){\circle*{2}}
\qbezier[1200](80,100)(79,106)(75,113)
\put(82,103){$\gamma$}\put(81,112){$\gamma$} \put(198,151){$z_2$}
\put(50,90){$a$} \put(350,90){$b$} \put(180,50){Figure 2}
\end{picture}

\end{document}